\documentclass[letterpaper, 10 pt, conference, onecolumn]{ieeeconf}

\IEEEoverridecommandlockouts
\usepackage{amsmath,amssymb,amsfonts}
\usepackage{algorithmic}
\usepackage{graphicx}
\usepackage{textcomp}
\usepackage{xcolor}
\usepackage{algorithm}
\usepackage{cite}
\usepackage{subcaption}

\allowdisplaybreaks

\def\BibTeX{{\rm B\kern-.05em{\sc i\kern-.025em b}\kern-.08em
    T\kern-.1667em\lower.7ex\hbox{E}\kern-.125emX}}

\newcommand{\cS}{\mathcal{S}}

\newcommand{\cA}{\mathcal{A}}

\newcommand{\cC}{\mathcal{C}}

\newcommand{\cU}{\mathcal{U}}
\newcommand{\cK}{\mathcal{K}}

\newcommand{\R}{\mathbb{R}}

\newcommand{\x}{x}

\newcommand{\z}{z}

\newcommand{\pr}{\epsilon}

\newcommand{\iter}{{t}}

\newcommand{\ud}{^}

\newcommand{\xt}{\x(\iter)}

\newcommand{\zt}{\z(\iter)}

\newcommand{\xtp}{\dot{\x}(t)}

\newcommand{\ztp}{\dot{\z}(\iter)}

\newcommand{\tz}{\tilde{\z}}

\newcommand{\tztp}{\dot{\tz}(\iter)}

\newcommand{\tzt}{\tz(\iter)}

\newcommand{\tzta}{\tz(\tau)}

\newcommand{\zeq}{\z\ud{\text{eq}}}

\newcommand{\norm}[1]{\left \|#1 \right \|}

\newcommand{\T}{^\top}

\newcommand{\map}[3]{#1: #2 \rightarrow #3}

\newcommand{\col}{\textsc{col}}

\newcommand{\lip}{L}

\newcommand\oprocendsymbol{\hbox{$\square$}}
\newcommand\oprocend{\relax\ifmmode\else\unskip\hfill\fi\oprocendsymbol}
\def\eqoprocend{\tag*{$\square$}}
\def\er/{Erd\H{o}s-R\'enyi}

\newcommand{\uu}{u}
\newcommand{\ut}{\uu(\iter)}

\newcommand{\dyn}{f}
\newcommand{\sdyn}{f} 
\newcommand{\fdyn}{g}

\newcommand{\rs}{_{\text{rs}}}
\newcommand{\bl}{_{\text{bl}}}

\newcommand{\nablahbound}{\lip_{\cbf}}
\newcommand{\fredbound}{\lip\rs}

\newcommand{\drs}{\sdyn\rs}
\newcommand{\dbl}{\fdyn\bl}

\def\PL/{Polyak-\L{}ojasiewicz}

\newcommand{\xddots}{%
	\raise 4pt \hbox {.}
	\mkern 6mu
	\raise 1pt \hbox {.}
	\mkern 6mu
	\raise -2pt \hbox {.}
}

\graphicspath{{figs/}}

\def\er/{Erd\H{o}s-R\'enyi}

\newcommand{\scale}{1}

\newcommand{\safeset}{\cS_x}
\newcommand{\safesetint}{\cS}

\newcommand{\cbf}{h}

\newcommand{\pol}{\pi} 
\newcommand{\lipp}{\lip_{\pol}}

\newcommand{\scl}{\tilde{\sdyn}_\pol}
\newcommand{\fcl}{\tilde{\fdyn}_\pol}
\newcommand{\zeqcl}{\zeq_{\pol}}

\newcommand{\out}{a}

\newcommand{\rob}{\eta}

\newcommand{\cCf}{\cC_{\text{\tiny full}}}

\def\citeDL/{\cite[Proposition~6.1.2]{bertsekas2015convex}}
 
\newtheorem{theorem}{Theorem}

\newtheorem{definition}{Definition}
\newtheorem{lemma}{Lemma}

\newtheorem{assumption}{Assumption}

\title{
  Safe Control of Feedback-Interconnected Systems via Singular Perturbations
  }
  \author{Stefano Di Gregorio, Guido Carnevale, Giuseppe Notarstefano
  \thanks{
    Work partially funded by FISA-2023-00210 project APACHE - CUP J53C25000520001, by the European Union - Next Generation EU - under the National Recovery and Resilience Plan (NRRP), Mission 4, Component 2, Investment 3.3; CUP J33C24001490009 and by IMA S.p.A.}
    \thanks{
      The authors are with the Department of Electrical, 
      Electronic and Information Engineering, 
      Alma Mater Studiorum Universit\`a di Bologna, 
      Bologna, Italy, {\tt\small name.lastname@unibo.it}}
      }

\begin{document}

\maketitle

\begin{abstract}
  Control Barrier Functions (CBFs) have emerged as a powerful tool in
  the design of safety-critical controllers for nonlinear systems.
  In modern applications, complex systems often involve the feedback
  interconnection of subsystems evolving at different timescales, e.g., two parts from different physical domains (such as the
  electrical and mechanical parts of robotic systems) or a physical
  plant and an (optimization or control) algorithm. In these
  scenarios, safety constraints often involve only a portion of the
  overall system.
  Inspired by singular perturbations for stability analysis, we
  develop a formal procedure to lift a safety certificate designed on
  a reduced-order model to the overall feedback-interconnected
  system.
  Specifically, we show that under a sufficient timescale separation
  between slow and fast dynamics, a composite CBF can be designed to
  certify the forward invariance of the safe set for the interconnected
  system.
  As a result, the online safety filter only needs to be solved
  for the lower-dimensional, reduced-order model.
  We numerically test the proposed approach on: (i) a robotic arm with
  joint motor dynamics, and (ii) a physical plant driven by an optimization
  algorithm.
\end{abstract}

\section{Introduction}

Autonomous systems operating in real world
scenarios require control algorithms capable of enforcing safety constraints.
This requirement has given rise to intense
research efforts in the literature, see the recent survey
\cite{hsu2023safety}.

Among existing methods, Control Barrier Functions (CBFs) have emerged as a
powerful tool for designing safe controllers (see the comprehensive
references~\cite{ames2019control,xiao2023safe,wang2025safe}). 
Their versatility has been shown across different domains, including
robotics~\cite{ferraguti2022safety}, autonomous
vehicles~\cite{alan2023control}, and human-robot
interaction~\cite{landi2019safety}.
A major limitation arises when the control input does not directly affect the safety constraint, leading to a high relative degree. 
To address this, High-Order CBFs (HOCBFs) \cite{xiao2021high} have been proposed to propagate safety via successive differentiations of the constraint.
Despite its effectiveness, this approach requires a uniform relative degree across the safe set.
This limitation has motivated a shift toward more flexible architectures, such as those based on input-output feedback linearization (see, e.g.,~\cite{cohen2024constructive}), rectified CBFs~\cite{ong2024rectified}, or the backstepping~\cite{freeman1993backstepping}.

In parallel, the integration of safety with other objectives has been explored. 
The combination of Control Lyapunov Functions (CLFs) and CBFs, allows for the simultaneous pursuit of safety and stability objectives. 
To this end, authors in \cite{ames2016control} unify them under a single optimization-based framework. 
This unification has been further generalized by the development of a universal formula for safe stabilization~\cite{ong2019universal}, providing an analytical construction for the control law. 
To handle practical uncertainties, these formulations have been extended to enforce safety and stability under input disturbances using input-to-state safe CBFs~\cite{kolathaya2018input}.
The conflict between reaching a target state and maintaining safety has been addressed in~\cite{matias2025hybrid}, where a hybrid CLF-CBF formulation is proposed to overcome intrinsic topological deadlocks.
Recently, an extension to model uncertainties~\cite{kamaldar2026composite} is derived via an adaptive approach obtained by combining a CBF, a CLF, and a parameter error term.

Beyond structural control design, a second major research direction focuses on the direct integration of CBFs into online dynamic optimization algorithms. 
In this framework, safety-critical control is reformulated as a constrained nonlinear program, where gradient flows are augmented with CBF-based inputs to ensure forward invariance of the feasible set \cite{allibhoy2023control}. 
This approach has been extended to the context of feedback optimization (see, e.g., \cite{colombino2019online, hauswirth2020timescale,carnevale2024nonconvex}). 
Specifically, authors of~\cite{delimpaltadakis2026safe} propose a safe feedback optimization method that enforces state constraints at all times.

However, the increasing complexity of modern architectures leads to nontrivial and/or unknown interconnected dynamics, for which designing CBFs is challenging or even impossible.
To mitigate this issue, there is a growing effort to develop
approaches synthesizing CBFs for reduced-order models, see the
recent survey~\cite{cohen2024safety}.
Several approaches rely on specific system structures. 
For instance, when the system is in strict-feedback form, safety can be propagated from a virtual controller to the actual input via backstepping~\cite{taylor2022safe}.
Alternative strategies include energy-based augmentations~\cite{singletary2021safety} or model-free tracking of safe reference velocities~\cite{molnar2021model}.
Closer to our methodology, the work~\cite{dall2026local} addresses safety in networked systems via a local implementation of the networked safety filter by resorting to the system structure typically used in singular perturbations.

Despite these advancements, a systematic framework to guarantee the
safety of generally interconnected systems, especially when dealing
with reduced-order models embedded within complex feedback loops,
remains an open challenge.

The main contribution of this paper is the establishment of a formal framework, inspired by singular perturbation (or timescale separation) theory for stability (see, e.g., the pioneering work~\cite{kokotovic1968singular} or the surveys~\cite{kokotovic1976singular,kokotovic1999singular,abdelgalil2023multi}), to guarantee safety of feedback-interconnected systems. 
Timescale separation arises naturally in modern applications in which there is interaction between physical domains evolving at different rates (e.g., the electrical and mechanical components of robotic systems) or between physical plants and optimization algorithms.
Specifically, we provide sufficient conditions to lift a CBF designed for a reduced-order model to the overall nonlinear, feedback-interconnected system.
Inspired by techniques from singular perturbations for stability, we show that if the neglected subsystem evolves on a sufficiently faster timescale than the reduced-order model, then a composite CBF can be constructed to certify the forward invariance of the safe set for the full interconnected system. 
Notably, our approach applies to broad classes of nonlinear dynamics without requiring a strict cascade structure or restrictive inversion-based feedback laws.
From an implementation perspective, the proposed framework allows safe controllers to be synthesized solely for the lower-dimensional reduced-order model. This significantly simplifies the safety filter design and reduces the online computational load, while still providing safety guarantees for the overall interconnection.

The paper is organized as follows. 
In Section~\ref{sec:preliminaries}, we review the definition of CBFs for continuous-time systems and the related safety guarantees.
In Section~\ref{sec:framework}, we introduce the considered framework and our goal.
In Section~\ref{sec:examples}, we present motivating examples for our framework.
In Section~\ref{sec:result}, we establish our main result. %
Finally, in Section~\ref{sec:simulations}, we test the proposed approach with numerical simulations.

\paragraph*{Notation}
$\R_{+}$ is the set $\{x \! \in \! \R \! \mid \! x\ge 0\}$.
We denote by $\col(v_1,\dots,v_N)$ the vertical concatenation of the vectors $v_1, \dots, v_N$.
A continuous function $\map{\alpha}{\R_+}{\R_+}$ is a class-$\cK_\infty$ function, denoted by $\alpha \in \cK_\infty$, if it is strictly increasing, $\alpha(0) = 0$, and $\lim_{r \to \infty} \alpha(r) = \infty$.

\section{Preliminaries on Safe Control}
\label{sec:preliminaries}

We start by recalling the theoretical foundations of CBFs for continuous-time nonlinear systems as
\begin{align}\label{eq:generic_system}
  \xtp = \dyn\left(\xt,\ut\right),
\end{align}
where $\xt \in \R^{n}$ is the system state at time $\iter \in \R_+$, $\ut \in \cU \subseteq \R^{m}$ is the input at time $\iter$, and $\dyn: \R^{n} \times \R^{m} \to \R^{n}$ describes the system dynamics.
The next definition introduces the notion of zero-superlevel set of a function.
\begin{definition}\label{def:zero_superlevel}
  Given $\cbf: \R^{n} \to \R$, the zero-superlevel set of $\cbf$ is defined as $\cC := \{\x \in \R^{n} \mid \cbf(\x) \geq 0\}$.
  \oprocend
\end{definition}
Then, let us 
introduce $\safeset \subseteq \R^n$ %
as the safe set of system~\eqref{eq:generic_system}.
Accordingly, we say that system~\eqref{eq:generic_system} is safe with respect to $\safeset$ if there exists a set $\cC \subseteq \safeset$ that is forward invariant for system~\eqref{eq:generic_system}.
A popular tool for constructively certifying the safety of system~\eqref{eq:generic_system} is given by CBFs \cite{ames2016control, ames2019control, cohen2024safety}.

For completeness, we report as follows a result to ensure safety via CBFs.
\begin{lemma}
  \label{lemma:cbf}
  \cite[Thm.~4]{cohen2024safety}
  Consider a continuously differentiable function $\cbf: \R^{n} \to \R$ with zero-superlevel set $\cC$. Assume that $\dyn$ is locally Lipschitz.
  If a locally Lipschitz feedback controller $\pol: \R^{n} \to \cU$ satisfies
  \begin{align*}
    \nabla \cbf(x)\T \dyn(\x,\pol(x)) \geq -\alpha(\cbf(\x)),
  \end{align*}
  for all $x \in \cC$, then $\cC$ is forward invariant for the closed-loop system
   \begin{align*}
    \xtp = \dyn(\xt,\pi(\xt)). \eqoprocend
  \end{align*}
\end{lemma}

\section{Singular Perturbations for Safety: Problem Statement}
\label{sec:framework}

\begin{figure}[H]
  \centering
  \includegraphics[scale=\scale]{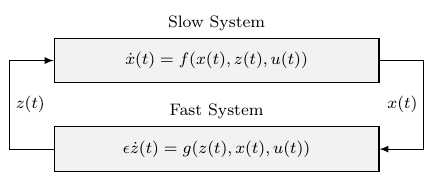}
  \caption{Block diagram representation of~\eqref{eq:interconnected_system_generic}.}
  \label{fig:scheme}
\end{figure}

In this section, we outline our procedure based on singular perturbations to lift reduced-order safety certificates in feedback interconnected systems.
Indeed, we focus on systems in the form
\begin{subequations}\label{eq:interconnected_system_generic}
  \begin{align}
    \xtp &= \sdyn(\xt,\zt,\ut)\label{eq:slow_system_generic}
    \\
    \pr\ztp &= \fdyn(\zt,\xt,\ut),\label{eq:fast_system_generic}
  \end{align}
\end{subequations}	
where $\xt \in \R^{n}$ and $\zt \in \R^{p}$ are the system states at
time $\iter \in \R_+$, $\ut \in \cU \subseteq \R^{m}$ is the input applied at time
$\iter$, $\sdyn: \R^n \times \R^p \times \R^m \to \R^n$ and
$\fdyn: \R^p \times \R^n \times \R^m \to \R^p$ are the vector fields
of the slow and fast subsystem, respectively, and $\pr > 0$ is a tunable parameter to
modulate the relative speed between these two subsystems.
First of all, in the next assumption, we impose Lipschitz continuity of $\sdyn$ and $\fdyn$.
\begin{assumption}\label{ass:lip}
  The functions $\sdyn$ and $\fdyn$ are Lipschitz continuous with
  parameters $\lip_{\sdyn}, \lip_{\fdyn} > 0$, respectively.
  \oprocend
\end{assumption}
By borrowing the typical nomenclature of timescale separation (or singular perturbations) for stability analysis (see, e.g., the recent survey~\cite{abdelgalil2023multi} or~\cite[Ch.~11]{khalil2002nonlinear}), we call $\xt$ and~\eqref{eq:slow_system_generic} \emph{slow} state and \emph{slow} system, respectively, $\zt$ and~\eqref{eq:fast_system_generic} are the \emph{fast} state and \emph{fast} system, respectively.
Fig.~\ref{fig:scheme} describes system~\eqref{eq:interconnected_system_generic} in terms of block diagrams.

\subsection{Reduced-Order Safety Certificates}

We are interested in designing safety certificates on reduced-order
models and then lift them to the overall feedback-interconnected
system. Purposely, we rely on the following three
requirements.

First, the safety specifications are defined solely in terms of the slow subsystem state $\x$ in \eqref{eq:slow_system_generic}.
In detail, we denote by
$\safeset \subseteq \R^{n}$ the safe set for the slow subsystem~\eqref{eq:slow_system_generic} and, accordingly, we define the overall safe set of
system~\eqref{eq:interconnected_system_generic} as
\begin{align}\label{eq:safesetint} 
  \safesetint := \{(x,z) \in \R^{n} \times \R^{p} \mid x \in \safeset\}.
\end{align}
Second, the fast subsystem~\eqref{eq:fast_system_generic} admits an \emph{equilibrium function} parametrized in $\x$ and $u$.
\begin{assumption}\label{ass:equilibrium}
  There exists a continuously differentiable function $\zeq: \safeset \times \cU \to \R^{p}$ such that
  \begin{align*}%
    0 &= \fdyn(\zeq(\x,u),\x,u).
  \end{align*}
  for all $\x \in \safeset$ and $u \in \cU$. \oprocend 
\end{assumption}
Third, we assume the availability of a robustly-safe controller
for a \emph{reduced-order model}
of~\eqref{eq:slow_system_generic}. Formally, we consider
system~\eqref{eq:interconnected_system_generic} on the manifold
$\{(x,z,u) \in \safeset \times \R^{p} \times \cU \mid z =
\zeq(x,u)\}$.
Namely, we introduce the reduced dynamics
$\drs: \R^n \times \R^m \to \R^n$ defined as
\begin{align}\label{eq:reduced}
  \drs(x,\uu) := \sdyn(x,\zeq(x,u),\uu).
\end{align}
For the reduced system we let the following assumption hold.

\begin{assumption}\label{ass:cbf}
  There exist a continuously differentiable function
  $\cbf: \R^{n} \to \R$ with a compact zero-superlevel set
  $\cC \subseteq \safeset$, a $\lipp$-Lipschitz continuous and
  differentiable function $\pol: \cC \to \cU$, a
  $\lip_{\alpha}$-Lipschitz continuous function
  $\alpha \in \cK_\infty$, and a safety margin $\rob > 0$ such
  that
  \begin{align}\label{eq:CBF}
    \nabla\cbf(x)\T \drs(x,\pol(x)) &\geq - \alpha(\cbf(x)) + \rob,
  \end{align}
  for all $x \in \cC$.\oprocend
\end{assumption}

We remark that Assumption~\ref{ass:cbf} provides a \emph{robust} safety certificate through the strictly positive safety margin $\rob>0$.
As shown later, this margin is instrumental in lifting the safety guarantee established for the reduced-order model~\eqref{eq:reduced} to the overall interconnected system~\eqref{eq:interconnected_system_generic}.
Indeed, the goal of this paper is to derive sufficient conditions under which the certificate in Assumption~\ref{ass:cbf} ensures forward invariance of a suitable set $\cCf \subseteq \safesetint$ (see the definition of $\safesetint$ in~\eqref{eq:safesetint}) for the overall interconnected system~\eqref{eq:interconnected_system_generic}. 
From an implementation standpoint, this allows us to enforce safety through a simplified filter designed solely for the reduced-order model, significantly reducing the online computational burden and eliminating the need for exact knowledge of the full feedback-interconnected dynamics.

\subsection{Boundary Layer System}

\begin{figure}[H]  
  \centering
  \includegraphics[scale=\scale]{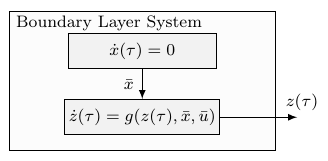}
  \caption{Block diagram representation of~\eqref{eq:bl} in original coordinates $(x,z)$.}
  \label{fig:scheme_bl}
\end{figure}

Here, following the nomenclature used in singular perturbations
  (timescale separation) for stability (see,
e.g.,~\cite{abdelgalil2023multi}), we characterize the behavior of
system~\eqref{eq:interconnected_system_generic} where
$\xt$ is fixed for all $\iter \in
\R_+$. %
We introduce an auxiliary system called \emph{boundary
  layer system} obtained (i) by arbitrarily \emph{freezing} the
state-input pair $(\bar{x},\bar{u}) \in \safeset \times \cU$ into the fast
dynamics~\eqref{eq:fast_system_generic}, (ii) by using the stretched
fast time variable $\tau := t/\pr$, and (iii) by considering the error
coordinates $\tzta := \z(\tau) - \zeq(\bar{x},\bar{u}) \in \R^{p}$ with respect
to the corresponding $\zeq(\bar{x},\bar{u})$.
Hence, the boundary layer system reads as
\begin{equation}\label{eq:bl}
	\frac{\partial \tzta}{\partial \tau} = \dbl(\tzta,\bar{x},\bar{u}),
\end{equation}
where $\tzt \in \R^p$ and
$\dbl: \R^p \times \R^{n} \times \R^m \to \R^p$ reads as
\begin{align}\label{eq:dbl_def}
  \dbl(\tz,\x,\uu) := \fdyn(\tz + \zeq(\x,\uu),\x,\uu).%
\end{align}
As depicted in Fig.~\ref{fig:scheme_bl}, the boundary-layer system~\eqref{eq:bl} represents~\eqref{eq:interconnected_system_generic} in the limit $\pr \to 0$, which is equivalent to freezing the slow state and the input at a given $(\bar{x}, \bar{u}) \in \safeset \times \cU$.
The next assumption ensures global exponential stability of the origin,
uniformly in $(\bar{x},\bar{u})$, for system~\eqref{eq:bl}.
\begin{assumption}\label{ass:bl}
  There exists a continuous function $U: \R^p \to \R$ such that
	\begin{subequations}\label{eq:U_generic}
		\begin{align}
			b_1 \norm{\tz}^2 \leq U(\tz) &\leq b_2\norm{\tz}^2
      \label{eq:U_first_bound_generic}
			\\
      \nabla U(\tz)\T \dbl(\tz,\bar{\x},\bar{\uu}) 
      &\leq 
      - b_3\norm{\tz}^2\label{eq:U_minus_generic}
			\\
      \norm{\nabla U(\tz)} &\leq b_4\norm{\tz},
      \label{eq:lipshcitz_hessian_u}
		\end{align}
	\end{subequations}
	for all $\tz \in \R^p$, $\bar{\x} \! \in \! \safeset$, $\bar{u} \in \cU$, and some $b_1, b_2, b_3, b_4 > 0$. \oprocend
\end{assumption}

\section{Motivating Examples Scenarios}\label{sec:examples}

Here, we provide two examples to motivate our setup.

\subsection{Robotic Manipulator with Actuator Dynamics}
\label{eq:motor_pendulum}

Consider a fully actuated 2-DoF planar robotic arm, where timescale separation arises naturally between mechanical and electrical domains.
Specifically, the mechanical dynamics of the system are described by the standard rigid-body equations
\begin{align}
  \label{eq:pendulum}
  \begin{bmatrix}
    \dot{q} \\
    \dot{\omega}
  \end{bmatrix} =
  \begin{bmatrix}
      \omega \\
      M(q)^{-1} \Big( K_I I - B \omega - C(q, \omega)\omega - G(q) \Big)
    \end{bmatrix},  
\end{align}
where $q = [q_1, q_2]^\top \in \mathbb{T}^2$ and $\omega = [\omega_1, \omega_2]^\top \in \R^2$ are angular positions and velocities, and $I \in \R^2$ represents motor currents. 
The terms $M$, $C$, $G$, $B$, and $K_I$ follow standard definitions in robotics.
Coupled with the mechanical plant, the electrical dynamics of the two motors is 
\begin{align}
  \label{eq:motor}
    L \dot{I} = V - RI - K_\omega \omega,
  \end{align}
where $V \in \R^2$ is the voltage applied to the motors, while $L$, $R$, $K_\omega \in \R^{2 \times 2}$ are positive definite diagonal matrices representing the inductance, resistance, and back-electromotive force constants, respectively.
Treating the inductance $L$ as the tunable parameter $\pr$ in~\eqref{eq:interconnected_system_generic}, we map this setup into our framework
with $x := \col(q, \omega)$, $\z := I$, $\uu := V$, and equilibrium function $I_\text{eq}: \mathbb{T}^2 \times \R^2 \times \R^2 \to \R^2$ defined as
\begin{align}
  I_\text{eq}(\col(q,\omega),V) := R^{-1} \left( V - K_\omega \omega \right).
\end{align}
Since $R \succ 0$, each equilibrium $I_\text{eq}([q,\omega]\T,V)$ is globally exponentially stable uniformly in $(q,\omega, V)$.
Hence, Assumptions~\ref{ass:equilibrium} and~\ref{ass:bl} hold true.
Assumption~\ref{ass:cbf} is verified by finding a safe controller on the reduced mechanical dynamics, namely 
\begin{align}\label{eq:motor_pendulum_red}
  \begin{bmatrix}
    \dot{q} 
    \\
    \dot{\omega}
  \end{bmatrix}
   =
  \begin{bmatrix}
    \omega \\
    M(q)^{-1} \Big( K_I R^{-1} V - B_{f} \omega - C(q, \omega)\omega - G(q) \Big)
  \end{bmatrix},
\end{align}
where $B_{f} := B + K_I R^{-1} K_\omega \in \R^{2 \times 2}$ is the effective damping matrix incorporating the back-electromotive force. %

\subsection{Suboptimal Safe Control}
\label{sec:dynamic_safe_control}

Second, we consider a setup in which the slow subsystem~\eqref{eq:slow_system_generic} is a given plant, while the fast one~\eqref{eq:fast_system_generic} is an optimization algorithm aimed at finding a safe control input that minimizes the distance from a desired but possibly unsafe control action. 
We formulate this scenario by introducing, for all $t \in \R_+$, an optimization problem of the form
\begin{subequations}\label{eq:problem}
  \begin{align}
    &\min_{u \in \cU} \quad \ell(u,\xt) 
    \\
    &\text{s.t. } \quad \nabla \cbf(\xt)\T \sdyn(\xt,u) \geq -\alpha(\cbf(\xt)) + \rob,  
  \end{align}
\end{subequations}
where $\ell: \R^m \times \R^n \to \R$ is the cost function
for a desired policy $\pi^{\text{des}}: \R^n \to \cU$, and $\cbf: \R^n \to \R$ is a CBF for $\sdyn$.
Denoting with $\cA: \R^{p} \times \R^{n} \to \R^{p}$ the algorithm flow used to address \eqref{eq:problem}, the interconnected system reads as 
\begin{subequations}\label{eq:plant_algo}
  \begin{align}
    \xtp &= \sdyn(\xt,\out(\zt))\label{eq:plant}
    \\
    \pr\ztp &= \cA(\zt,\xt),\label{eq:algo}
  \end{align}
\end{subequations}
where $\pr > 0$ tunes the algorithm speed~\eqref{eq:algo}, and $\out: \R^p \to \cU$ extracts the current solution estimate from $\zt$.
In this setup, Assumption~\ref{ass:cbf} is satisfied when problem~\eqref{eq:problem} admits a solution function $\uu_\star: \cC \to \cU$, that is continuously differentiable in the set $\cC$.
Assumption~\ref{ass:equilibrium} is satisfied provided that $\cA$ admit an equilibrium function $\zeq: \cC \to \R^p$ such that $\out(\zeq(x)) = \uu_\star(x)$ for all $x \in \cC$.
Assumption~\ref{ass:bl} is achieved with exponential stability of $\zeq(x)$ for system $\ztp = \cA(\zt,x)$ uniformly in $\x \in \cC$.
Finally, Assumption~\ref{ass:lip} follows from the Lipschitz continuity of $\sdyn$ and $\cA$.

\section{Singular Perturbations for Safety: Main Result}
\label{sec:result}

In this section, we provide a safety certificate for the interconnected system~\eqref{eq:interconnected_system_generic} based on the reduced-order safety $\cbf$ and the stability properties of the boundary layer system~\eqref{eq:bl} (cf. Assumption~\ref{ass:bl}), respectively.
To this end, we introduce the composite CBF function $V: \cC \times \R^{p} \to \R$ defined as
\begin{align}
  V(x,z) := \cbf(x) - U(z - \zeq(x,\pol(x))),\label{eq:V}
\end{align}
and we denote as $\cCf := \{(x,z) \in \R^{n} \times \R^{p} \mid V(x,z) \geq 0\}$ its zero-superlevel set.
\begin{theorem}\label{th:SP}
	Consider system~\eqref{eq:interconnected_system_generic} and let Assumptions~\ref{ass:lip},~\ref{ass:equilibrium},~\ref{ass:cbf}, and~\ref{ass:bl} hold.
  Then, there exists $\bar{\pr} > 0$ such that, for all $\pr \in (0,\bar{\pr})$, the set $\cCf$ is forward invariant for the closed-loop system
    \begin{align*}
      \xtp &= \sdyn(\xt,\zt,\pol(\xt))
      \\
      \pr\ztp &= \fdyn(\zt,\xt,\pol(\xt)).\eqoprocend
    \end{align*}
\end{theorem}
The proof of Theorem~\ref{th:SP} is provided in Appendix~\ref{sec:proof}.

\iftrue
To illustrate the theoretical guarantees provided by Theorem \ref{th:SP}, we consider the 2D singularly perturbed system
\begin{subequations}\label{eq:toy}
  \begin{align}
    \dot{x}(t) &= z(t) - x(t),\label{eq:slow_toy} 
    \\
    \epsilon \dot{z}(t) &= u(t) + x(t) - z(t),\label{eq:fast_toy}
  \end{align}
\end{subequations}
with $x(t)$, $z(t)$, $\ut \in \R$.
The fast subsystem equilibrium function is $\zeq(x,u) = u + x$. 
We consider the safe set 
\begin{align*} 
  \safeset = \{x \in \R \mid \beta^2 - x^2 \geq 0\},
\end{align*}
for some $\beta > 0$, and the CBF $\cbf: \R \to \R$ defined as
\begin{align*}
  \cbf(x) := \beta^2 - x^2.
\end{align*}
Namely, in this case, $\cC \equiv \safeset$.
To achieve this safety goal while tracking a nominal policy $\pol_{\text{des}}(x) = 2 \sin(x)$, the control input $u(t) = \pol(x)$ is computed via a CBF-QP based on the reduced-order dynamics, enforcing a strict safety margin $\rob>0$ to the constraint.

In Fig.~\ref{fig:cbf_comparison}, we provide graphical representations that help visualize the concepts highlighted below regarding the meaning of Theorem~\ref{th:SP}.

\begin{figure}[H]
  \centering
  \includegraphics[width=0.48\textwidth]{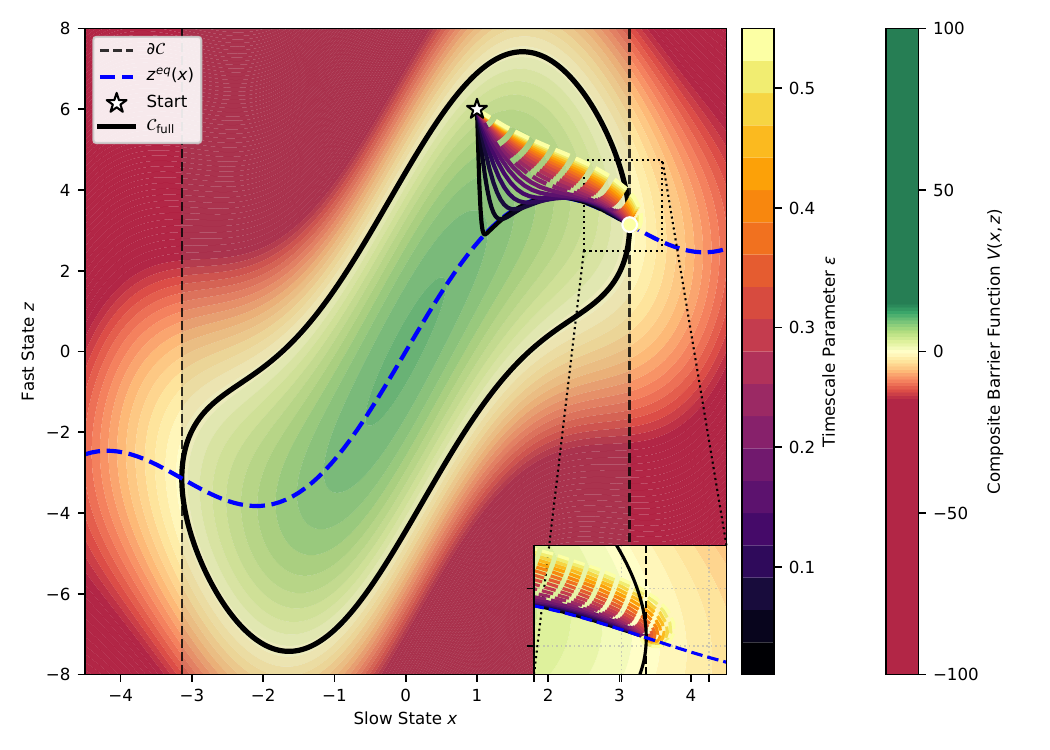}
  \caption{Graphical representation of the safety guarantees provided by Theorem~\ref{th:SP} for a 2D example. The phase plane displays the level sets of the composite CBF $V(x,z)$, the composite safe set boundary $\cCf$ (solid black line), the ideal reduced-order safe set boundary (dashed vertical lines), and the equilibrium manifold $\zeq(x,\pol(x))$ (dashed blue curve). Trajectories are depicted for varying values of the timescale parameter $\pr$. Safe trajectories are shown with solid lines, while unsafe trajectories are represented by dashed lines.
  }
  \label{fig:cbf_comparison}
\end{figure}

\paragraph{Invariant Set Lifting}
  
  As shown in Fig.~\ref{fig:cbf_comparison}, the composite set $\cCf$ (enclosed by the solid black line) is strictly contained within the bounds $[-\beta, \beta]$ of the ``ideal'' reduced-order safe set (dashed vertical lines). 
  Theorem~\ref{th:SP} ensures that starting within $\cCf$, the interconnected system remains safe for sufficiently small $\pr$. 
  This introduces a restriction: $\cbf(x)$ must be large enough to compensate for the fast state deviation $z - \zeq(x, \pol(x))$, causing the effective safe region for $x$ to shrink as $z$ moves away from the equilibrium manifold (dashed blue curve).

\paragraph{Timescale Separation and Safety}

  Theorem \ref{th:SP} captures the required timescale separation. 
  As $\pr \to 0$, the fast state instantly converges to $\zeq(x,\pol(x))$, preserving safety. 
  For $\pr \neq 0$, the transient mismatch acts as a perturbation, and the strict safety margin $\rob$ serves as a robustness buffer. 
  When $\pr < \bar{\pr}$, this perturbation is absorbed, maintaining safety.
  If $\pr \geq \bar{\pr}$, the deviation may lead to constraint violations.

\fi

\section{Numerical Simulations}\label{sec:simulations}

To validate the theoretical results of Theorem~\ref{th:SP}, we perform numerical simulations on (i) the fully actuated 2-DoF robotic arm described in Section~\ref{eq:motor_pendulum}, and (ii) the suboptimal safe control scenario described in Section~\ref{sec:dynamic_safe_control}.
\subsection{Robotic Arm with Joint Motors}

Consistently with our theoretical framework, we design the safety certificate strictly for the reduced-order mechanical dynamics~\eqref{eq:motor_pendulum_red}. 
The physical parameters used in the simulations are provided in Table \ref{tab:robot_parameters} (both the links and the motors are assumed to be identical to one another).
\begin{table}[H]
\centering
\caption{Manipulator and Actuator Parameters}
\label{tab:robot_parameters}
\begin{tabular}{|l|c|c|c|c|c|c|c|c|}
\hline
\textbf{Symbol} & $m$ & $l$ & $J$ & $b$ & $R$ &$K_I$ & $K_\omega$ \\ \hline
\textbf{Value}  & 0.5 & 1.0 & 0.167 & 0.5 & 0.6 & 0.4 & 0.4        \\ \hline
\end{tabular}
\end{table}
As discussed in Section~\ref{sec:examples}, the inductance $L$ is treated as a tunable parameter $\pr$ to modulate the timescale separation between the mechanical and electrical dynamics.
We define a set of $N_c=9$ safety margins $\bar{h}_j(q) \geq 0$ to encompass both the joint limits and the obstacle avoidance constraints.
In detail, 
given 
$q_{\min} = [0^\circ, -80^\circ]^\top$ and $q_{\max} = [185^\circ, 80^\circ]^\top$, the first $4$ functions $\bar{h}_j(q)$ are defined as
\begin{align*}
  \bar{h}_1(q) &:= q_1 - q_{\min,1}
  ,\quad
  \bar{h}_2(q) := q_{\max,1} - q_1 
  \\
  \bar{h}_3(q) &:= q_2 - q_{\min,2} 
  ,\quad
  \bar{h}_4(q) := q_{\max,2} - q_2.
\end{align*}
Further, we consider $5$ circular obstacles with radius $r = 0.3$m centered at the workspace coordinates $c_1 = [1.1, -1.7]^\top$, $c_2 = [2.2, 0.0]^\top$, $c_3 = [1.9, 0.6]^\top$, $c_4 = [1.3, 1.7]^\top$, and $c_5 = [-0.4, 2.1]^\top$.
To avoid them, we define %
\begin{align*}
  \bar{h}_{4+k}(q) := \norm{p(q) - c_k}^2 - r^2, \quad k \in \left\{1, \dots, 5\right\},
\end{align*}
where $p(q) \in \R^2$ is the position of the end-effector.
Since these constraints have a relative degree of two with respect to the input, we construct a High-Order CBF (see, e.g., \cite{xiao2021high}) for each $j\in\{1,\dots,N_c\}$, namely, given $\gamma > 0$, we define
\begin{align}
  h_j(q, \omega) := \nabla \bar{h}_j(q)\T \omega + \gamma \bar{h}_j(q).
\end{align}
To simplify the filter design, all constraints are aggregated into a single CBF $h(q, \omega)$ using the Log-Sum-Exp formulation (see, e.g., \cite{safari2025safe}), namely
\begin{align}
  h(q, \omega) := -\tfrac{1}{\rho} \ln \left(\sum\textstyle_{j=1}^{N_c} e^{-\rho h_j(q, \omega)}\right),
\end{align}
where $\rho > 0$ modulates the approximation sharpness. 
Then, we choose a controller that guarantees the safety certificate
\begin{align}
  \nabla h(q, \omega)\T f_{\text{rs}}(q, \omega) + \gamma h(q, \omega) - \rob \geq 0,
\end{align}
where $f_{\text{rs}}$ is the vector field of the reduced-order system~\eqref{eq:motor_pendulum_red}. 

\begin{figure*}[htpb]
  \centering
  \begin{subfigure}[b]{0.32\textwidth}
    \centering
    \includegraphics[width=\textwidth]{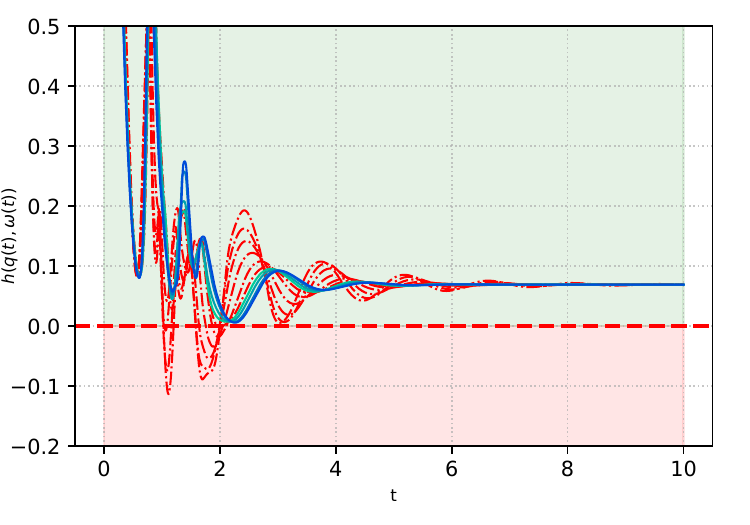}
    \caption{}
    \label{fig:h_dp}
  \end{subfigure}
  \hfill
  \begin{subfigure}[b]{0.32\textwidth}
    \centering
    \includegraphics[width=\textwidth]{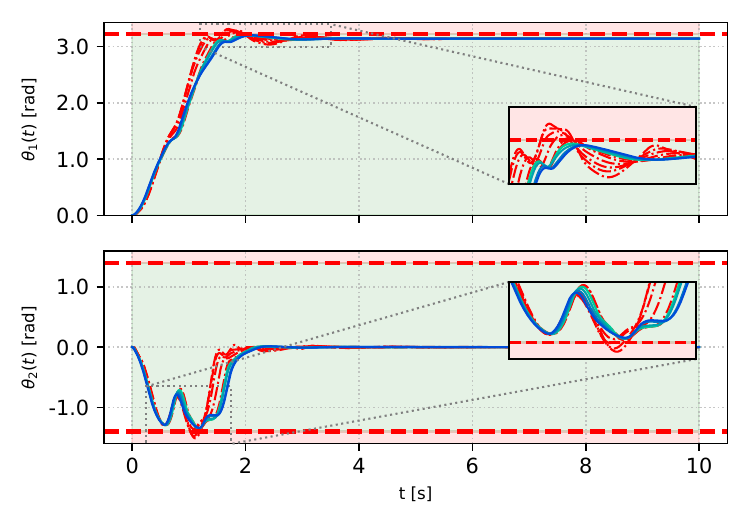}
    \caption{}
    \label{fig:theta_dp}
    \end{subfigure}
    \hfill
    \begin{subfigure}[b]{0.32\textwidth}
      \centering
      \includegraphics[width=\textwidth]{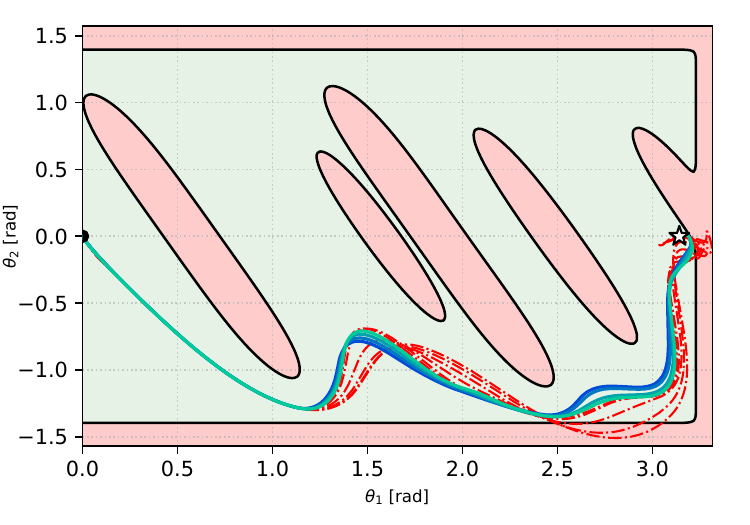}
      \caption{}
      \label{fig:cspace_dp}
      \end{subfigure}
  \caption{Results for the 2-DoF robotic arm with joint motors. Time evolution of $\cbf(q, \omega)$ (a) and of $q = [\theta_1, \theta_2]^\top$ (b). System trajectories in the joint space (c), where red areas indicate unsafe regions. Blue-to-cyan solid lines indicate safe trajectories (corresponding to increasing $\pr < \bar{\pr}$), whereas the dash-dotted red lines denote safety violations ($\pr \geq \bar{\pr}$).}
  \label{fig:dp_history}
\end{figure*}

We numerically implement a smooth safety filter, following the line of \cite{cohen2023characterizing}. The system, initialized at rest ($q(0) = [0^\circ, 0^\circ]^\top$, $\omega(0) = [0, 0]^\top$), is driven toward the target configuration $q_\text{target} = [180^\circ, 0^\circ]^\top$ by a nominal PD controller with gravity compensation. 
The parameters of the safety filter are set as $\gamma = 10$, $\rob = 0.1$, and $\rho = 50$.
Simulations are performed for $10$ linearly spaced values of $\pr \in \left[0.001, 0.035\right]$.
Results are shown in Fig. \ref{fig:dp_history}. As theoretically predicted by Theorem \ref{th:SP}, the interconnected system remains safe for sufficiently small $\pr$. 
Conversely, as $\pr$ exceeds a critical threshold, constraint violations occur.

\subsection{Suboptimal Safe Control}

\begin{figure}[H]
  \centering
  \includegraphics[width=0.48\textwidth]{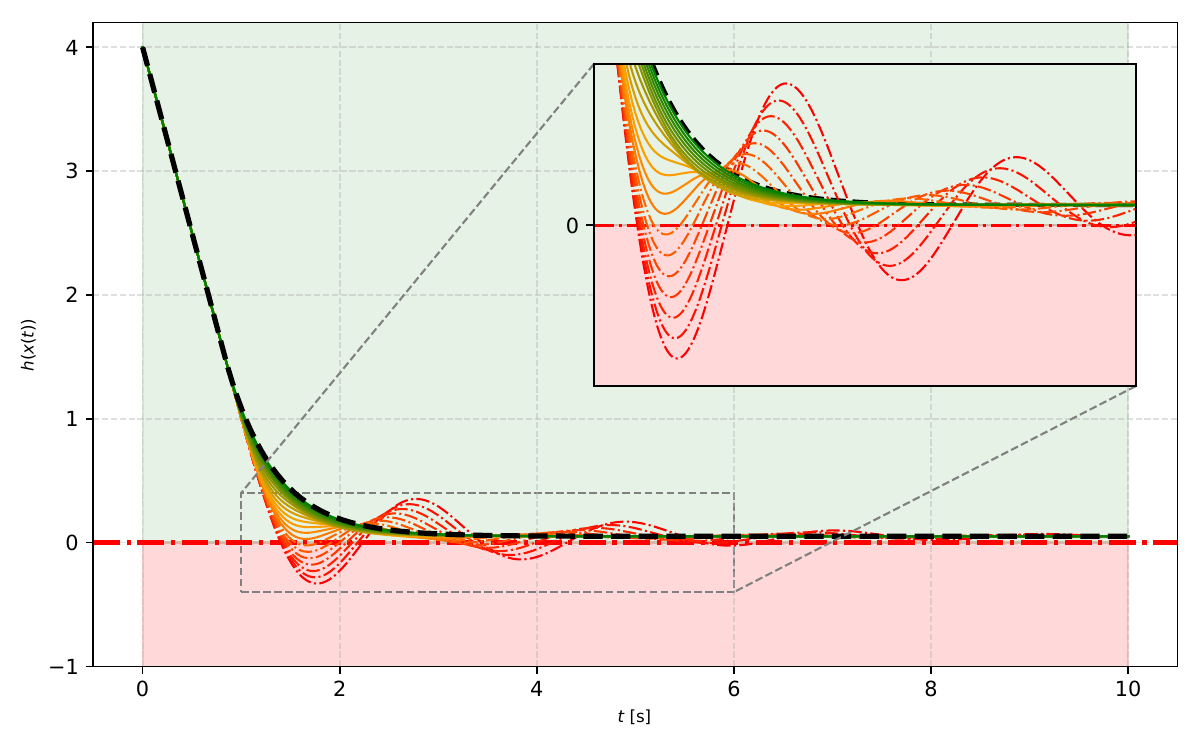}
  \caption{Time evolution of the CBF $\cbf(x)$ for varying values of $\pr$. The green-to-red color gradient indicates increasing values of $\pr$, with the ideal CBF-QP solution represented by the dashed black line. Unsafe trajectories are dash-dotted.}
  \label{fig:h_evolution}
\end{figure}

We consider a system $\dot{x} = u$ and
the CBFs $\cbf_1(x) := x$ and $\cbf_2(x) := 10 - x$. As in the previous example, we aggregate them into a single CBF $\cbf(x) := -\frac{1}{\rho} \ln (e^{-\rho \cbf_1(x)} + e^{-\rho \cbf_2(x)})$. 
We apply an unsafe nominal policy $\pi^{\text{des}} = -3$. To dynamically implement the safety filter, we define the algorithm $\cA(z, x)$ in \eqref{eq:algo} as the primal-dual gradient flow
\begin{align*}
    \cA(z, x) := \begin{bmatrix}
    -\nabla_u \ell(z, x) - \lambda \nabla_u g(z, x) 
    \\
    \Pi_{\R_+}(g(z, x))
    \end{bmatrix},
\end{align*}
where $z := [u, \lambda]^\top$ comprises the control input and Lagrange multiplier, $\ell(u, x) = (u - \pi^{\text{des}})^2$ is the cost function, $g(u, x) := \nabla \cbf(x) u + \alpha(\cbf(x)) - \rob$ encodes the robust safety constraint, and $\Pi_{\R_+}$ is the projection onto $\R_+$. 
The system is initialized at $x(0) = 4$ and $z(0) = [\pi^{\text{des}}, 0]^\top$. 
The parameters of the safety filter are set as $\gamma = 2$, $\rho=50$, and $\rob = 0.1$.
Simulations span $10$ seconds for $20$ linearly spaced values of $\pr \in \left[0.01, 0.3\right]$.
We benchmark our suboptimal algorithm against an exact CBF-QP solution, computed by solving \eqref{eq:problem} at each time step. 
As shown in Fig. \ref{fig:h_evolution}, the CBF-QP strictly maintains $\cbf(x(t)) \geq 0$.
Conversely, the suboptimal approach preserves safety only for sufficiently small values of $\pr$. 

\section{Conclusions}

In this paper, we developed a singular perturbation framework for safety. %
Specifically, we considered feedback interconnections between two subsystems, referred to as the slow and fast subsystems, respectively.
We then provided sufficient conditions to guarantee safety of the overall interconnected system, based on safety properties of the reduced subsystem and stability properties of the fast one.
We illustrated the applicability of our results in two settings.
The first is the safe control of a robotic arm with joint motors.
The second is a suboptimal safe control scenario in which a physical plant interacts with an algorithm that iteratively solves an optimization problem encoding safety constraints.

\appendix

\subsection{Proof of Theorem~\ref{th:SP}} 
\label{sec:proof}

The idea of the proof is to combine the reduced-system safety properties encoded in $\cbf$ (cf. Assumption~\ref{ass:cbf}) with the boundary-layer system~\eqref{eq:bl} stability properties encoded in $U$ (cf. Assumption~\ref{ass:bl}) in order to guarantee that $V(x,z) = \cbf(x) - U\bigl(z - \zeq(x,\pol(x))\bigr)$ ensures safety for the overall interconnected system~\eqref{eq:interconnected_system_generic}.
To this end, the proof exploits the parameter $\pr$ to arbitrarily slow down the evolution of the slow state $\xt$.
Intuitively, this allows one to make arbitrarily small the approximation errors arising from the fact that the properties of $\cbf$ and $U$ are established only for the limit cases detailed in Assumptions~\ref{ass:cbf} and~\ref{ass:bl}, rather than for the original interconnection~\eqref{eq:interconnected_system_generic}.

First of all, since they will be useful later, let us introduce the function $\zeqcl: \cC \to \R^{p}$ and the constants $\nablahbound, \fredbound, \lip_{\zeqcl}  > 0$ defined as
\begin{subequations}
  \begin{align}
    \zeqcl(x) &:= \zeq(x,\pol(x))
    \label{eq:zeqcl}
    \\
    \fredbound &:= \sup_{x \in \cC}  
    \norm{\sdyn(\x,\zeqcl(\x),\pol(\x))} 
    \label{eq:bound_phi}
    \\
    \nablahbound &:= \sup_{x \in \cC}  
    \norm{\nabla \cbf(x)}
    \label{eq:nabla_W_alpha_W}
    \\
    \lip_{\zeqcl} &:= \sup_{x \in \cC}\norm{\nabla \zeqcl(x)}.\label{eq:zeqcl_lip}
\end{align}
\end{subequations}
We highlight that the above constants exist and are finite since $\cC$ is compact (cf. Assumption~\ref{ass:cbf}) $\sdyn$, $\zeq$, $\nabla\zeq$, $\nabla \cbf$, and $\pol$ are continuous (see Assumptions~\ref{ass:lip},~\ref{ass:equilibrium}, and~\ref{ass:cbf}, respectively).
Then, by Lipschitz continuity of $\sdyn$ and $\fdyn$ (cf. Assumption~\ref{ass:lip}) and $\pol$ (cf. Assumption~\ref{ass:cbf}), we note that the solution to system~\eqref{eq:interconnected_system_generic} exists and is unique (see, e.g.,~\cite[Thm.~3.2]{khalil2002nonlinear}).
Let us rewrite system~\eqref{eq:interconnected_system_generic} in the error coordinate $\tzt := \zt - \zeq(\xt,\pol(\xt))$, thus obtaining 
\begin{subequations}\label{eq:interconnected_system_generic_err}
  \begin{align}
    \xtp &= \scl(\xt,\tzt)
    \label{eq:slow_system_generic_err}
    \\
    \tztp &= \frac{1}{\pr}\fcl(\tzt,\xt) -\nabla\zeqcl(\xt)\xtp
    ,\label{eq:fast_system_generic_err}
  \end{align}
\end{subequations}
in which we introduce $\scl: \cC \times \R^{p} \to \R^{n}$ and $\fcl: \R^p \times \cC \to \R^{p}$ defined as 
\begin{subequations}
  \begin{align}
   \scl(x,\tz) &:= \sdyn(x,\tz + \zeqcl(x),\pol(x))
    \label{eq:scl}
    \\
    \fcl(\tz,x) &:= \fdyn(\tz + \zeqcl(x),x,\pol(x)).
    \label{eq:fcl}
  \end{align}
\end{subequations}
Now, let us consider a pair $(x,\tz)$ such that $\cbf(x) - U(\tz) \ge 0$ and recall that this necessarily implies $x \in \safeset$ since $\cC \subseteq \safeset$ (cf. Assumption~\ref{ass:cbf}) and $U(\tz) \ge 0$ for all $\tz \in \R^{p}$ (see~\eqref{eq:U_first_bound_generic} in Assumption~\ref{ass:bl}).
Then, we add and subtract the term $\nabla\cbf(\x)\T\drs(\x,\pol(\x)) = \nabla\cbf(\x)\T\scl(\x,0)$ (see the definitions of $\drs$~\eqref{eq:reduced} and $\scl$~\eqref{eq:scl}) to the Lie derivative of $\cbf$ along the trajectories of subsystem~\eqref{eq:slow_system_generic_err}, thus obtaining
\begin{align}
  \nabla\cbf(\x)\T\scl(\x,\tz)
  &= 
  \nabla\cbf(\x)\T\drs(\x,\pol(\x))
  + \nabla\cbf(\x)\T(\scl(\x,\tz)-\scl(\x,0))
  \notag\\
  &\stackrel{(a)}{\geq}
  -\alpha(\cbf(\x)) + \rob
  - \norm{\nabla\cbf(\x)}\norm{\scl(\x,\tz)-\scl(\x,0)}
  \notag\\
  &\stackrel{(b)}{\geq}
  -\alpha(\cbf(\x)) + \rob
  - \nablahbound\lip_{\sdyn} \norm{\tz},
  \label{eq:deltaW}
\end{align}
where in $(a)$ we exploit~\eqref{eq:CBF} (cf. Assumption~\ref{ass:cbf}) since $\x \in \cC$, and we use the Cauchy-Schwarz inequality, while in $(b)$ we use the Lipschitz continuity of $\sdyn$ (cf. Assumption~\ref{ass:lip}) and the bound~\eqref{eq:nabla_W_alpha_W}.
Now, we evaluate the Lie derivative of the Lyapunov function $U$ characterized in Assumption~\ref{ass:bl} along the trajectories of~\eqref{eq:fast_system_generic_err} and we use~\eqref{eq:U_minus_generic} to write
\begin{align}
  \frac{1}{\pr}\nabla U(\tz)\T\fcl(\tz,\x) - \nabla U(\tz)\T\nabla\zeqcl(\x)\xtp
  &\leq 
  - \frac{b_3}{\pr}\norm{\tz}^2 - \nabla U(\tz)\T\nabla\zeqcl(\x)\xtp
  \notag\\        
  &\stackrel{(a)}{\leq}
  - \frac{b_3}{\pr}\norm{\tz}^2 + b_4\lip_{\zeqcl}\norm{\tz}\norm{\scl(\x,\tz)}
  \notag\\        
  &\stackrel{(b)}{\leq}
  - \frac{b_3}{\pr}\norm{\tz}^2 + b_4\lip_{\zeqcl}\norm{\tz}\norm{\drs(\x,\pol(\x))}
  +b_4\lip_{\zeqcl}\norm{\tz}\norm{\scl(\x,\tz) - \scl(\x,0)}\notag\\        
  &\stackrel{(c)}{\leq}
  - \frac{b_3}{\pr}\norm{\tz}^2 + b_4\lip_{\zeqcl}\fredbound\norm{\tz}
  +b_4\lip_{\zeqcl}\lip_{\sdyn}\norm{\tz}^2
  \label{eq:DeltaU}
\end{align} 
where in $(a)$ we use the Cauchy-Schwarz inequality, the bounds~\eqref{eq:lipshcitz_hessian_u} (cf. Assumption~\ref{ass:bl}) and~\eqref{eq:zeqcl_lip}, and~\eqref{eq:slow_system_generic_err}, in $(b)$ we add and subtract the term $\drs(\x,\pol(\x)) = \scl(x,0)$ (see the definitions of $\drs$~\eqref{eq:reduced} and $\scl$~\eqref{eq:scl}) in the last norm, and use the triangle inequality, while in $(c)$ we use the Lipschitz continuity of $\sdyn$ (cf. Assumption~\ref{ass:lip}) and the bound~\eqref{eq:bound_phi}. %
Now, let us adapt the definition of $V$ (cf.~\eqref{eq:V}) to the error coordinate $\tz$, namely, let us introduce $\tilde{V}: \cC \times \R^p \to \R$ defined as
\begin{align}
  \tilde{V}(\x,\tz) := V(x,\tz+\zeq(x)) = \cbf(\x) - U(\tz).\label{eq:tilde_V}
\end{align}
Then, the Lie derivative of $\tilde{V}$ along the trajectories of the interconnected system~\eqref{eq:interconnected_system_generic_err} reads as
\begin{align}
  \nabla\tilde{V}(\x,\tz)\T
  \begin{bmatrix}
    \scl(\x,\tz) \\
    \frac{1}{\pr}\fcl(\tz,\x) + \nabla\zeqcl(\x)\scl(\x,\tz)  
  \end{bmatrix}
  & 
  =\nabla \cbf(\x)\T \scl(\x,\tz) 
  - \frac{1}{\pr}\nabla U(\tz)\T\fcl(\tz,\x) -\nabla U(\tz)\T\nabla\zeqcl(\x)\scl(\x,\tz)   
  \notag\\
  &\stackrel{(a)}{\geq} 
  -\alpha(\cbf(\x)) + \rob + \left(\frac{b_3}{\pr} - b_4\lip_{\zeqcl}\lip_{\sdyn}\right)\norm{\tz}^2
  - (\nablahbound\lip_{\sdyn} + b_4\lip_{\zeqcl}\fredbound)\norm{\tz}
  \label{eq:deltaV_generic}
\end{align}
where in $(a)$ we substitute the expressions found in \eqref{eq:deltaW} and \eqref{eq:DeltaU}.
To simplify notation, let us introduce the constants $K_1 := \nablahbound\lip_{\sdyn} + b_4\lip_{\zeqcl}\fredbound$ and $K_2 := b_4\lip_{\zeqcl}\lip_{\sdyn}$, thus resulting in the following bound for the Lie derivative of $\tilde{V}$:
\begin{align}
  &\nabla\tilde{V}(\x,\tz)\T
  \begin{bmatrix}
    \scl(\x,\tz) \\
    \frac{1}{\pr}\fcl(\tz,\x) + \nabla\zeqcl(\x)\scl(\x,\tz)  
  \end{bmatrix} 
  \geq -\alpha(\cbf(\x)) + \rob + \left(\frac{b_3}{\pr} - K_2\right)\norm{\tz}^2 - K_1\norm{\tz}.\label{eq:delta_tilde_V}
\end{align}
Now, let us arbitrarily choose $\nu \in (0,1)$ and consider the decomposition 
\begin{align*}
  \frac{b_3}{\pr}\norm{\tz}^2 = \nu\frac{b_3}{\pr}\norm{\tz}^2 + (1-\nu)\frac{b_3}{\pr}\norm{\tz}^2.
\end{align*}
Then, by imposing the condition 
\begin{align*} 
  \nu\frac{b_3}{\pr} - K_2 > 0 \iff \pr < \bar{\pr}_1 := \frac{\nu b_3}{K_2}, 
\end{align*}
we observe that the terms $-K_1\norm{\tz} + \left(\nu\frac{b_3}{\pr} - K_2\right)\norm{\tz}^2$ form a parabola of $\norm{\tz}$ opening upwards  and.
Thus, we can bound
this parabola through its minimum value, namely 
\begin{align*} 
  -K_1\norm{\tz} + \left(\nu\frac{b_3}{\pr} - K_2\right)\norm{\tz}^2 \geq -\frac{K_1^2}{4\left(\nu\frac{b_3}{\pr} - K_2\right)},
\end{align*} 
which allow us to further bound the right-hand side of~\eqref{eq:delta_tilde_V} according to
\begin{align}
  \nabla\tilde{V}(\x,\tz)\T
  \begin{bmatrix}
    \scl(\x,\tz) \\
    \frac{1}{\pr}\fcl(\tz,\x) + \nabla\zeqcl(\x)\scl(\x,\tz)  
  \end{bmatrix}
    &
  \geq -\alpha(\cbf(\x)) + (1-\nu)\frac{b_3}{\pr}\norm{\tz}^2
  + \rob - \frac{K_1^2}{4\left(\nu\frac{b_3}{\pr} - K_2\right)}
  \notag\\
  &\stackrel{(a)}{\geq} 
  -\alpha(\cbf(\x)) + (1-\nu)\frac{b_3}{\pr}\norm{\tz}^2,  
  \label{eq:lie_V_generic}
\end{align}
where in $(a)$ we neglect the term $\rob - \frac{K_1^2}{4\left(\nu\frac{b_3}{\pr} - K_2\right)}$ by making it positive through the choice 
\begin{align*} 
  \pr < \bar{\pr}_2 := \min\left\{\bar{\pr}_1, \tfrac{4 \rob \nu b_3}{K_1^2 + 4 \rob K_2}\right\}.
\end{align*}
Now, by $\lip_\alpha$-Lipschitz continuity of $\alpha$ (cf. Assumption~\ref{ass:cbf}), it holds
\begin{align}\label{eq:lip_alpha_generic}
  |\alpha(\cbf(\x)) - \alpha(\cbf(\x) - U(\tz))| \leq \lip_{\alpha}|U(\tz)|
\end{align}
Since $U(\tz) \ge 0$ (see~\eqref{eq:U_first_bound_generic} in Assumption~\ref{ass:bl}), and $\cbf(\x) - U(\tz) \ge 0$ by hypothesis, the bound~\eqref{eq:lip_alpha_generic} leads to
\begin{align}
  \alpha(\cbf(\x)) &\leq \alpha(\cbf(\x) - U(\tz)) + \lip_{\alpha}U(\tz)
  \notag\\
  &\stackrel{(a)}{\leq} \alpha(\tilde{V}(\x,\tz)) + \lip_{\alpha}b_2\norm{\tz}^2,
  \label{eq:lipschitz_alpha}
\end{align}
where in $(a)$ we use the definition of $\tilde{V}$ (cf.~\eqref{eq:tilde_V}) and the bound~\eqref{eq:U_first_bound_generic}.
Then, by using~\eqref{eq:lipschitz_alpha}, we can further bound the right-hand side of~\eqref{eq:lie_V_generic} as 
\begin{align}
  \nabla\tilde{V}(\x,\tz)\T
  \begin{bmatrix}
    \scl(\x,\tz) \\
    \frac{1}{\pr}\fcl(\tz,\x) + \nabla\zeqcl(\x)\scl(\x,\tz)  
  \end{bmatrix}
  &\geq
  -\alpha(\tilde{V}(\x,\tz)) - \lip_{\alpha}b_2\norm{\tz}^2 + (1-\nu)\frac{b_3}{\pr}\norm{\tz}^2,\label{eq:last} 
\end{align}
By setting $\bar{\pr} := \min\{\bar{\pr}_2, \frac{(1-\nu)b_3}{\lip_{\alpha}b_2}\}$ and using $\pr < \bar{\pr}$, we neglect that last two terms in the right-hand side of~\eqref{eq:last}, thus obtaining
\begin{align}
  &\nabla\tilde{V}(\x,\tz)\T
  \begin{bmatrix}
    \scl(\x,\tz) \\
    \frac{1}{\pr}\fcl(\tz,\x) + \nabla\zeqcl(\x)\scl(\x,\tz)  
  \end{bmatrix}
  \geq
  -\alpha(\tilde{V}(\x,\tz)).
\end{align}
The proof follows by using the last bound to invoke Lemma~\ref{lemma:cbf}.

\end{document}